\documentclass[12pt]{amsart}
\usepackage{amssymb}
\usepackage{verbatim}
\usepackage[usenames]{color}
\usepackage{hyperref}

\newtheorem{thm}{Theorem}
\newtheorem{prop}[thm]{Proposition} 
\newtheorem{la}[thm]{Lemma}
\newtheorem{cor}[thm]{Corollary}

\theoremstyle{definition}
\newtheorem{df}[thm]{Definition}

\newtheorem{notat}[thm]{Notation}
\newtheorem{qu}[thm]{Question}

\theoremstyle{remark}
 
\newtheorem{rmk}[thm]{Remark}

\newenvironment{ls}{\begin{itemize}}{\end{itemize}}
\newenvironment{lsnum}{\begin{enumerate}}{\end{enumerate}}
\newenvironment{pf}{\begin{proof}}{\end{proof}}
\newcommand{\ger}[1]{\ensuremath{\mathfrak {#1}}}
\newcommand{\scr}[1]{\ensuremath{\mathcal {#1}}}
\newcommand{\bld}[1]{\ensuremath{\mathbf {#1}}}
\newcommand{\bbb}[1]{\ensuremath{\mathbb {#1}}}

\newcommand{\emp}{\varnothing}
\renewcommand{\phi}{\varphi}

\newcommand{\sq}[1]{\ensuremath{\langle#1\rangle}}

\newcommand{\restr}{\mathop{\upharpoonright}}

\newcommand{\notarrow}{\kern .42em\not\kern -.42em\longrightarrow}
\renewcommand{\th}{\ensuremath{{}^{\text{th}}}}

\newcommand{\sm}{{\text-}\!\!\sum}
\newcommand{\fr}{\ensuremath{\mathfrak{Fr}}}
\newcommand{\pr}{\text{-prod}\,}
\newcommand{\prn}{\pr\ger N}

\newcommand{\noprint}[1]{\relax}

\title{Partitions and Conservativity}
\author{Andreas Blass}
\address{Mathematics Department\\
University of Michigan\\
Ann Arbor, MI 48109--1043, U.S.A.}
\email{ablass@umich.edu}

\begin{document}

\begin{abstract}
  We study the partition properties enjoyed by the ``next best thing
  to a P-point'' ultrafilters introduced recently in joint work with
  Dobrinen and Raghavan.  That work established some finite-exponent
  partition relations, and we now analyze the connections between
  these relations for different exponents and the notion of
  conservativity introduced much earlier by Phillips.  In addition, we
  establish some infinite-exponent partition relations for these
  ultrafilters and also for sums of non-isomorphic selective
  ultrafilters indexed by selective ultrafilters.
\end{abstract}

\maketitle

\section{Introduction}

This paper is, in part, a sequel to an earlier joint paper with
Natasha Dobrinen and Dilip Raghavan \cite{nbt} in which we studied the
ultrafilters created by a certain forcing construction.  These
ultrafilters were shown to have numerous combinatorial properties,
including some weak partition properties and a classification of
arbitrary functions modulo the ultrafilter.  

Part of the present paper is an analysis of the implications between
these combinatorial properties in general, i.e., independently from
the connection with the generic ultrafilters of \cite{nbt}.  This
analysis involves also a notion of conservativity that was introduced
by Phillips \cite{phillips}, was studied further in \cite{ee-ce-rf,
  model-th}, and is connected with the model-theoretic notion of
stability.

The other part of the present paper concerns infinitary partition properties
of the generic ultrafilters in \cite{nbt}.  These partition properties
lead to so-called ``complete combinatorics'' for the forcing notion
used in \cite{nbt}.  A different approach to these partition
properties was developed independently by Dobrinen \cite{dobrinen}.

\section{Fr\'echet-Squared Genericity}

In this section, we review the results from \cite{nbt} that will play
a major role in the present paper, and we take the opportunity to also
fix some notation and terminology.

We deal with filters on countably infinite sets, and, in the absence
of a contrary statement, filters are assumed to extend the cofinite
filter.  Sometimes we formulate definitions and results in terms of
filters on the set $\omega$ of natural numbers, but these definitions
and results can be transferred to other countably infinite sets via
arbitrary bijections.

\begin{df}              \label{p+sel:df}
  An ultrafilter \scr U on $\omega$ is a \emph{P-point} if, for every
  function $f$ on $\omega$, there is some $A\in\scr U$ such that the
  restriction $f\restr A$ is finite-to-one or constant. \scr U is
  \emph{selective} if, for every function $f$ on $\omega$, there is
  some $A\in\scr U$ such that the restriction $f\restr A$ is
  one-to-one or constant.
\end{df}

Clearly, every selective ultrafilter is a P-point.  The existence of
selective ultrafilters as well as the existence of non-selective
P-points can be proved if the continuum hypothesis is assumed, but it
is consistent with ZFC that there are no P-points.  

We shall need a well-known eqivalent characterization of P-points.  An
ultrafilter \scr U is a P-point if, given any countably many sets
$A_n\in\scr U$, there is a set $B\in\scr U$ almost included in all of
them, i.e., $B-A_n$ is finite for all $n$.  The proof of the
equivalence is based on making an $f$ in the definition have value $n$
on $A_n-A_{n-1}$ (where $A_{-1}$ means $\emp$).

Selective ultrafilters enjoy remarkable partition properties and are
therefore sometimes called Ramsey ultrafilters.  Specifically, we have
the following result, in which the part about partitions of
$[\omega]^n$ is due to Kunen (published in Booth's thesis
\cite{booth}) and the parts about $[\omega]^\omega$ are due to Mathias
\cite{happy}. 

\begin{prop}            \label{sel-ramsey}
Every selective ultrafilter \scr U on $\omega$ has the following
partition properties.
\begin{lsnum}
\item For all $n,k\in\omega$, we have $\omega\to(\scr U)^n_k$, which
  means that, whenever the set $[\omega]^n$ of $n$-element subsets of
  $\omega$ is partitioned into $k$ pieces, then there is a set
  $H\in\scr U$ such that $[H]^n$ is included in one piece.
\item $\omega\overset{\text{analytic}}\longrightarrow(\scr
  U)^\omega_2$, 
  which means that, whenever the set $[\omega]^\omega$ of infinite
  subsets of $\omega$ is partitioned into an analytic piece and its
  complement, then there is a set $H\in\scr U$ such that
  $[H]^\omega$ is included in one piece.
\item Suppose the universe is obtained from some ground model by
  L\'evy-collapsing to $\omega$ all cardinals below a Mahlo cardinal
  of the ground model.  Then $\omega\overset{HOD\bbb
    R}\longrightarrow(\scr U)^\omega_2$, where $HOD\bbb R$ means
  hereditarily ordinal definable from reals.
\end{lsnum}
\end{prop}

In part~(2) of this proposition, ``analytic'' refers to the topology
of $[\omega]^\omega$ as a subspace of the power set $\scr P(\omega)$,
which is, in turn, identified with $2^\omega$ topologized as a product
of discrete two-point spaces.  Thus, two infinite subsets of $\omega$
are near each other in $[\omega]^\omega$ just in case they have a long
common initial segment.

The partition relation $\omega\to(\scr U)^n_k$ for any particular $n$
and $k\geq2$ easily implies the same relation for the same $n$ and all
$k$; it also implies the same relation for all smaller $n$.  Also,
$\omega\to(\scr U)^2_2$ easily implies selectivity, because, given a
function $f$ on $\omega$, we can partition the set $[\omega]^2$ of
pairs $\{a,b\}$ according to whether $f(a)=f(b)$.  The proposition
thus implies the additional fact that the partition relation
$\omega\to(\scr U)^n_k$ for any one value of $n$ implies the same for
larger values of $n$, and this is not such an easy result.

The most natural way to produce a nonprincipal ultrafilter \scr U on
$\omega$ by forcing is to take, as forcing conditions, the infinite
subsets of $\omega$, ordered by inclusion (so smaller sets are
stronger conditions).  The intended meaning of a condition $A$ is that
it forces $A\in\scr U$.  This ordering is not separative; the
separative quotient identifies two infinite sets if and only if their
symmetric difference is finite, so it amounts to the Boolean
algebra\footnote{When we use Boolean algebras as notions of forcing,
  we always mean the algebras with their zero elements removed.} $\scr
P(\omega)/\text{fin}$.  This separative quotient is countably closed,
so the forcing adds no new reals.  That makes it easy to check, by a
density argument, that the generic object added by the forcing is an
ultrafilter and indeed a selective ultrafilter.  

Part~(3) of Proposition~\ref{sel-ramsey} has an important consequence
concerning this $\scr P(\omega)/\text{fin}$ forcing.  First, let us
weaken that part of the proposition as follows.  In the L\'evy-Mahlo
model (i.e., the model obtained by L\'evy collapsing all cardinals
below some Mahlo cardinal to $\omega$), if $\scr D\subseteq
[\omega]^\omega$ is in $HOD\bbb R$, then either (a)~there is an infinite
$H\subseteq\omega$ such that no infinite subset of $H$ is in \scr D,
or (b)~there is an $H\in\scr U\cap\scr D$.  Note that, in both
alternatives, we have weakened the conclusion in
Proposition~\ref{sel-ramsey}, part~(3).  In (a), we allow $H$ to be
any infinite set, not necessarily in \scr U.  In (b), we only require
$H$ itself, not all its infinite subsets, to be in \scr D.  This
weakened form of the partition relation can be succinctly restated as
follows. 

\begin{cor}             \label{comcom1}
In the L\'evy-Mahlo model, every selective ultrafilter is $\scr
P(\omega)/\text{fin}$-generic over $HOD\bbb R$.
\end{cor}

Because of this corollary, one says that selectivity is \emph{complete
  combinatorics} for $\scr P(\omega)/\text{fin}$ forcing.  Any generic
ultrafilter for this forcing has the combinatorial property of
selectivity, and there is no genuinely\footnote{The word ``genuinely''
  here is intended to exclude information like the partition
  properties in Proposition~\ref{sel-ramsey}, which intuitively look
  like more information than selectivity but in fact follow from
  selectivity.} additional combinatorial information that follows from
$\scr P(\omega)/\text{fin}$-genericity, because, at least in the
L\'evy-Mahlo model, selectivity by itself already ensures genericity
over a large inner model.

\begin{rmk}             \label{forcing-equiv} 
The forcing notion $\scr
  P(\omega)/\text{fin}$ is equivalent, for forcing, to the set of
  countably generated filters on $\omega$ (extending the cofinite
  filter, as usual), partially ordered by reverse inclusion.  Indeed,
  we can associate to any infinite subset $A$ of $\omega$ the filter
  $\{X\in[\omega]^\omega:A\subseteq^*X\}$, where $\subseteq^*$ means
  inclusion modulo finite sets.  This gives a dense embedding of $\scr
  P(\omega)/\text{fin}$ into the partially ordered set of countably
  generated filters.  

  The analog of Proposition~\ref{sel-ramsey} part~(3) for this partial
  ordering of filters is false.  Specifically, the image of our
  embedding of $\scr P(\omega)/\text{fin}$ is, as mentioned above,
  dense, but so is the complement of this image.  That is, every
  countably generated filter can be extended to one that does not
  contain a set $A$ that is included modulo finite in all other
  elements of the extended filter.

  Curiously, though, despite the failure of
  Proposition~\ref{sel-ramsey} part~(3) for the filter poset, the
  weakened version used in proving Corollary~\ref{comcom1} is true.
  One way to see this is to observe that this weakened version is
  equivalent to the corollary, which, being only about forcing, is
  clearly unchanged when we replace the forcing notion $\scr
  P(\omega)/\text{fin}$ by an equivalent one.
\end{rmk}

One of the objectives of the paper \cite{nbt} was to obtain, for a
suitable ultrafilter that is not a P-point, results analogous to those
for selective ultrafilters described above.  To this end, we
introduced what seems to be the simplest notion of forcing that
adjoins a non-P-point ultrafilter on $\omega$, and we studied the
properties of these ultrafilters in considerable detail.  In order to
describe this forcing and the resulting ultrafilters, it is convenient
to first introduce some notation, which will also be useful later in
other contexts.

\begin{df}              \label{sums:df}
For a subset $A$ of $\omega^2$, we define its \emph{(vertical)
  sections} to be the sets 
\[
A(x)=\{y\in\omega:\sq{x,y}\in A\}.
\]
For filters \scr U and $\scr V_n$ on $\omega$, we define the \scr
U-indexed \emph{sum} of the $\scr V_n$'s to be
\[
\scr U\sm_n\scr V_n=\{A\subseteq\omega^2:
\{n\in\omega:A(n)\in\scr V_n\}\in\scr U\}.
\]
That is, a subset of the plane is large with respect to this sum if
and only if almost all (with respect to \scr U) of its sections are
large (with respect to the appropriate $\scr V_n)$.  It is easy to
verify that this sum is always a filter, and that it is an ultrafilter
if \scr U and all of the $\scr V_n$'s are ultrafilters.

If all of the $\scr V_n$'s are the same filter \scr V, then we write
$\scr U\otimes\scr V$ for the sum. If, furthermore, $\scr V=\scr U$,
then we use the notation $\scr U^{\otimes 2}$.  

Finally, we write \fr\ for the filter of cofinite subsets of $\omega$,
because it is often called the Fr\'echet filter.
\end{df}

The forcing notion that formed the main subject of \cite{nbt} is the
Boolean algebra\footnote{I prefer to think of quotients of Boolean
  algebras as being determined by filters rather than by ideals.  So I
  use the filter notation $\fr^{\otimes 2}$ here, even though I used
  the standard notation $\scr P(\omega)/\text{fin}$ earlier rather
  than my preferred $\scr P(\omega)/\fr$.} $\scr
P(\omega^2)/\fr^{\otimes2}$.  It is shown in \cite{nbt} that this is a
countably closed forcing notion and therefore does not add reals.  As
a consequence, the generic object that it adjoins is an ultrafilter
\scr G on $\omega^2$ extending the filter $\fr^{\otimes 2}$.  This
implies that, for every $A\in\scr G$, there are infinitely many
$n\in\omega$ such that $A_n$ is infinite; indeed, this property of $A$
is exactly what it means for $A$ to intersect every set in
$\fr^{\otimes2}$ or, equivalently, not to be in the ideal dual to
$\fr^{\otimes2}$.

A frequently useful notion of forcing equivalent to $\scr
P(\omega^2)/\fr^{\otimes2}$ can be obtained as follows.  First,
replace the equivalence classes that constitute the quotient algebra
$\scr P(\omega^2)/\fr^{\otimes2}$ with all their representatives.
That is, form the poset of all subsets of $\omega^2$ that meet every
set in $\fr^{\otimes2}$, i.e., that have infinitely many infinite
sections.  This poset is not separative; its separative quotient is
$\scr P(\omega^2)/\fr^{\otimes2}$.  Second, pass to the sub-poset
consisting of those elements that have no finite sections.  This is a
dense sub-poset, and in fact its separative quotient is still $\scr
P(\omega^2)/\fr^{\otimes2}$, because, when we remove all finite
sections from a set, we do not alter its equivalence class modulo
$\fr^{\otimes2}$. Third, pass to the even smaller poset consisting of
those conditions on which the second projection $\omega^2\to\omega$ is
one-to-one.  This sub-poset is still dense in our forcing, because,
given any condition $A$, we can thin out all its sections so as to be
disjoint from each other yet still infinite.  Fourth, restrict to
those conditions that lie entirely above the diagonal, i.e., that
consist only of pairs $(x,y)$ with $x<y$.  This is again a dense
subset, and it makes no difference in the separative quotient since
the part of $\omega^2$ that lies above the diagonal is a set in
$\fr^{\otimes2}$. Finally, for technical convenience, consider only
those conditions $A$ such that we never have the $x$-coordinate of one
point in $A$ equal to the $y$-coordinate of another point in $A$.
Conditions of this sort are easily seen to be dense in our poset, so
the final result is still forcing equivalent to the original $\scr
P(\omega^2)/\fr^{\otimes2}$. Summarizing this construction, we have
that $\scr P(\omega^2)/\fr^{\otimes2}$ is equivalent, as a forcing
notion, to the poset defined as follows.

\begin{df}              \label{p:df}
  \bbb P is the set of those subsets $A$ of $\omega^2$ that satisfy
  \begin{lsnum}
    \item $A$ has infinitely many infinite sections and no nonempty
      finite sections.
\item The sections of $A$ are pairwise disjoint.
\item All elements \sq{x,y} of $A$ have $x<y$.
\item For any \sq{x,y} and \sq{x',y'} in $A$, we have $x\neq y'$.
  \end{lsnum}
\end{df}
We sometimes identify the above-diagonal subset of $\omega^2$,
$\{\sq{x,y}:x<y\}$, with the set $[\omega]^2$  of two-element subsets
$\{x<y\}$ of $\omega$.  Thus, the forcing conditions in \bbb P can be
viewed as subsets of $[\omega]^2$, and the generic \scr G can be
viewed as an ultrafilter on $[\omega]^2$.

An immediate consequence of the ``infinitely many infinite sections''
property of sets in \scr G is that the projection
$\pi_1:\omega^2\to\omega$ to the first factor is neither finite-to-one
nor constant on any set in \scr G.  Therefore, the generic filter is
not a P-point.  

Here is a list of additional properties of $\scr
P(\omega^2)/\fr^{\otimes2}$-generic ultrafilters\footnote{Genericity
  here is over $V$, which I think of as the universe of all sets, so
  generic objects are in Boolean extensions. The results remain
  correct under any of the other customary ways to view forcing, for
  example taking $V$ to be a countable transitive model so that
  generic objects exist in the ordinary, two-valued universe of sets.}
proved in \cite[Section~3]{nbt}; the references are to the
propositions, theorems, and corollaries of that paper.  After the
list, I shall provide definitions for the concepts used in the list.

\begin{prop}            \label{nbt-prop}
If \scr G is $\scr P(\omega^2)/\fr^{\otimes2}$-generic over $V$, then
the following statements hold in $V[\scr G]$.
\begin{lsnum}
\item The ultrafilter $\pi_1(\scr G)$ is $\scr
  P(\omega)/\text{fin}$-generic over $V$ and therefore
  selective. (Proposition~30) 
\item The projection $\pi_2:\omega^2\to\omega$ to the second factor is
  one-to-one on a set in \scr G. (Corollary~32)
\item For any function $f$ on $\omega^2$, there is a set $A\in\scr G$
  such that $f\restr A$ is one of 
  \begin{ls}
    \item a constant function,
\item $\pi_1$ followed by a one-to-one function, and
\item a one-to-one function. (Corollary~33)
  \end{ls}
\item \scr G is a weak P-point.  (Theorem~36)
\item For any partition of $[\omega^2]^n$ into finitely many pieces,
  there is a set $H\in\scr G$ such that, for any $n$-type $\tau$, all
  $n$-element subsets of $[H]^n$ that realize $\tau$ lie in the same
  class of the partition.  Therefore, \scr G is $(n,T(n))$-weakly
  Ramsey, where $T(n)$ is the number of $n$-types.  (Theorem~31)
\end{lsnum}
\end{prop}

For the sake of completeness, we explain here the notation $\pi_1(\scr
G)$ and the terminology ``weak P-point'', ``$n$-type'', ``realize'',
and ``weakly Ramsey'' used in items~1, 4, and 5 of this proposition.
More information about these can be found in \cite{nbt}; the notion of
weak P-point comes from Kunen's paper \cite{kunen}, where he proved
(in ZFC) that weak P-points exist.

For any function $f:X\to Y$ and any filter \scr U on $X$, its image on
$Y$ is defined as 
\[
f(\scr U)=\{A\subseteq Y:f^{-1}(A)\in\scr U\}.
\]
This is always a filter, except that it may fail to satisfy our
convention that filters on $\omega$ must contain all cofinite sets.
If \scr U is an ultrafilter, then so is $f(\scr U)$, which extends the
cofinite filter as long as $f$ is not constant on any set in \scr U.

\begin{df}              \label{weak-p:df}
  An ultrafilter \scr U on a countably infinite set $S$ is a
  \emph{weak P-point} if, given any countable set of (nonprincipal)
  ultrafilters $\scr W_n$ on $S$, all distinct from \scr U, we have a
  set $A$ that is in \scr U but in none of the $\scr W_n$.
\end{df}

In terms of the topology of $\beta S-S$, this means that \scr U is not
in the closure of any countable set of ultrafilters distinct from \scr
U.  It is easy to verify that, as the terminology implies, all
P-points are weak P-points.  Kunen showed in \cite{kunen} that there
always exist weak P-points that are not P-points.

A sum $\scr U\sm_n\scr V_n$ is never a weak P-point, as can be seen by
taking the $\scr W_n$ in Definition~\ref{weak-p:df} to be copies of
the $\scr V_n$ on the vertical sections of $\omega^2$.  That is, let
$i_n:\omega\to\omega^2$ be the map $y\mapsto\sq{n,y}$ and set $\scr
W_n=i_n(\scr V_n)$.

\begin{df} \label{weak-ramsey:df} Let $n$ and $t$ be natural numbers,
  let $S$ be a countably infinite set, and let \scr U be a
  (nonprincipal) ultrafilter on $S$. Then \scr U is
  $(n,t)$-\emph{weakly Ramsey} if, whenever $[S]^n$ is partitioned
  into finitely many pieces, there is a set $H\in\scr U$ such that
  $[H]^n$ meets at most $t$ of the pieces.
\end{df}

If $t=1$, this is the partition property in part~(1) of
Proposition~\ref{sel-ramsey}.  As $t$ increases, the $(n,t)$-weak
Ramsey property gets weaker.  As $n$ increases, the property gets
stronger.  

Another common notation for this property is $S\to[\scr U]^n_{t+1}$.
The reason for the subscript $t+1$ is that this partition relation is
equivalent to saying that, whenever $[S]^n$ is partitioned into $t+1$
pieces, there is an $H\in\scr U$ such that $[H]^n$ is disjoint from at
least one piece.  Our definition, with an arbitrary finite number of
pieces follows easily from this version with just $t+1$ pieces, by
induction on the number of pieces.

Note that, in arrow notations like $S\to[\scr U]^n_{t+1}$, the square
brackets are used to indicate the weak form of homogeneity, merely
missing a piece, whereas round parentheses indicate the strong form,
meeting only one piece.  This notation is so common that one often
speaks of square-bracket partition relations.

To explain $n$-types, it is useful to begin by considering an
arbitrary element $A$ of the notion of forcing \bbb P from
Definition~\ref{p:df}.  Let us agree to write any $n$-element subset
of $A$ as $\{\sq{a_1,b_1},\dots,\sq{a_n,b_n}\}$ with
$b_1<\dots<b_n$. Recall that clause~(2) in Definition~\ref{p:df}
ensures that the $b_i$'s are all distinct, so we are merely adopting
the convention to list the $n$ pairs \sq{a_i,b_i} in the order of
increasing second components.  The $n$-type realized by $A$ will be
defined as all the information about the relative ordering of the
$a_i$'s and $b_i$'s, with no information about their actual values.
More precisely, we define types and realization as follows.

\begin{df}              \label{type:df}
An $n$-\emph{type} is a linear pre-order of the set of $2n$ formal
symbols $x_1,\dots,x_n,y_1,\dots,y_n$ such that 
\begin{ls}
\item $y_1<\dots<y_n$,
\item each $x_i$ precedes the corresponding $y_i$, and
\item if two distinct symbols are equivalent in the pre-order, then
  both of them are $x$'s.
\end{ls}
\end{df}
Recall that a pre-order on a set is a reflexive, transitive, binary
relation $\leq$ on that set; that it is linear if every two elements
of the set are ordered one way or the other; that two elements are
called equivalent if each is $\leq$ the other; that identifying
equivalent elements leads to a partial order (linear if the preorder
was linear) on the quotient set; and that $<$ means ``$\leq$ and not
$\geq$.''

In \cite[Definition~2.9]{nbt}, we used a different formulation of the
notion of $n$-type, namely a list of the $x_i$'s and $y_i$ with $=$ or
$<$ between each consecutive pair, subject to requirements
corresponding to the clauses of the present definition.  We pointed
out, in \cite[Remark~2.11]{nbt}, that this list form of $n$-types is
equivalent to the pre-order form adopted here.  My main reason for now
preferring the pre-order version is that is generalizes more naturally
to the case of infinite sets in place of 
$\{a_1,\dots,a_n,b_1,\dots,b_n\}$.  Nevertheless, the list form is
also convenient, for example in the following definition.

\begin{df}              \label{realize:df}
The $n$-type \emph{realized} by an $n$-element subset
$\{\sq{a_1,b_1},\dots,\sq{a_n,b_n}\}$ of an element of \bbb P is the
pre-order of $\{x_1,\dots,x_n,y_1,\dots,y_n\}$ whose list form becomes
true when the $x_i$'s and $y_i$'s are interpreted as denoting the
corresponding $a_i$'s and $b_i$. 
\end{df}

It is easy to verify, in the light of Definition~\ref{p:df} and our
convention that $b_1<\dots<b_n$, that every $n$-element subset of an
element of \bbb P realizes a unique type.  Furthermore, every
$A\in\bbb P$ has $n$-element subsets realizing all of the $n$-types.
This result is essentially Proposition~2.14 of \cite{nbt}, though it
is stated and proved in somewhat greater generality there.

\begin{notat}
  $T(n)$ denotes the number of $n$-types.
\end{notat}

The preceding discussion shows that the generic ultrafilter cannot be
$(n,T(n)-1)$-weakly Ramsey.  Indeed, the same goes for any ultrafilter
having a basis of sets from \bbb P and, as \cite[Corollary~2.16]{nbt}
shows, for any ultrafilter that is not a P-point.  Thus, part~5 of
Proposition~\ref{nbt-prop} says that \scr G has the strongest
weak-Ramsey properties that are possible for a non-P-point.

It will be convenient to have shorter name for the property in part~(3)
of Proposition~\ref{nbt-prop}. 

\begin{df}
  An ultrafilter \scr U on $\omega^2$ has the \emph{three functions
    property} if every function on $\omega^2$ is, when restricted to
  some set in \scr U, either one-to-one or constant or the composition
  of $\pi_1$ followed by a one-to-one function.
\end{df}

The motivation for this terminology is that there are, up to
restriction to sets in \scr U and post-composition with one-to-one
functions, just three functions on $\omega$, namely the identity,
$\pi_1$, and any constant function.  By analogy, selectivity could be
called the two-functions property.

\section{Selective-Indexed Sums of Selective Ultrafilters}

In this section, we describe another family of ultrafilters on
$\omega^2$ enjoying many but not all of the properties of the $\scr
P(\omega^2)/\fr^{\otimes2}$-generic (or equivalently \bbb P-generic)
ultrafilters discussed in the preceding section.  These ultrafilters
are the sums $\scr U\sm_n\scr V_n$, as in Definition~\ref{sums:df}, of
selective ultrafilters, in the special case that the summands are
pairwise non-isomorphic and the indexing ultrafilter \scr U is also
selective.  To avoid having to repeatedly use the long phrase
``\bld{s}elective-\bld{i}ndexed \bld{s}um of
\bld{n}on-\bld{i}somorphic \bld{s}electives'', we introduce the
following acronym.

\begin{df}
A \emph{sisnis} ultrafilter is an ultrafilter on $\omega^2$ of the
form 
\[
\scr U\sm_n\scr V_n
\]
where \scr U and all of the $\scr V_n$'s are selective ultrafilters
and, for $m\neq n$, there is no permutation $f$ of $\omega$ with
$f(\scr V_m)=\scr V_n$.
\end{df}
The content of the definition would be unchanged if we required only
\scr U-almost all of the $\scr V_n$'s to be selective and pairwise
non-isomorphic.  This is because a sum $\scr U\sm_n\scr V_n$ is
unchanged if we change the $\scr V_n$'s arbitrarily for a set of $n$'s
whose complement is in \scr U.  For the same reason, we can assume
that \scr U is not isomorphic to any of the $\scr V_n$'s, since this
can be arranged by altering at most one $\scr V_n$.  The content would
also be unchanged if we allowed $f$ to be an arbitrary function
$\omega\to\omega$ rather than a permutation.  This is because, by
selectivity, any $f:\omega\to\omega$ is $\scr V_m$-almost everywhere
equal to either a constant function or a one-to-one function.  A
constant $f$ cannot map $\scr V_m$ to a non-principal ultrafilter such
as $\scr V_n$, and a one-to-one map would be $\scr V_m$-almost
everywhere equal to a permutation of $\omega$.

The next proposition summarizes information from \cite[Section~2]{nbt}
about sisnis ultrafilters; see Lemmas~2.2 and 2.3, Proposition~2.4,
and Theorem~2.17 of \cite{nbt}.

\begin{prop} \label{sisnis-prop} Every sisnis ultrafilter $\scr W=\scr
  U\sm_n\scr V_n$ on $\omega^2$ has the following properties.
\begin{lsnum}
\item The ultrafilter $\pi_1(\scr W)$ is the selective ultrafilter
  \scr U.
\item The projection $\pi_2:\omega^2\to\omega$ to the second factor is
  one-to-one on a set in \scr W. 
\item For any function $f$ on $\omega^2$, there is a set $A\in\scr W$
  such that $f\restr A$ is one of 
  \begin{ls}
    \item a constant function,
\item $\pi_1$ followed by a one-to-one function, and
\item a one-to-one function. 
  \end{ls}
\item \scr W is not a weak P-point.
\item For any partition of $[\omega^2]^n$ into finitely many pieces,
  there is a set $H\in\scr W$ such that, for any $n$-type $\tau$, all
  $n$-element subsets of $[H]^n$ that realize $\tau$ lie in the same
  class of the partition.  Therefore, \scr W is $(n,T(n))$-weakly
  Ramsey.
\end{lsnum}
\end{prop}

Do not be lulled by the apparent similarity between this proposition
and Proposition~\ref{nbt-prop}.  Although the other clauses in the two
propositions match, clause~(4) is entirely different.  \bbb P-generic
ultrafilters are weak P-points but sisnis ultrafilters are
not. Indeed, as pointed out earlier, sums of ultrafilters are never
weak P-points.

Of course, it follows that sisnis ultrafilters \scr W are not
P-points.  It is easy to see this directly, because the first
projection $\pi_1$ is neither finite-to-one nor constant on any set in
\scr W.  

We need one additional property of sisnis ultrafilters, the analog of
a trivial property of \bbb P-generic ultrafilters.

\begin{la}
Every sisnis ultrafilter has a basis consisting of sets in \bbb P.
\end{la}

\begin{pf}
  Consider an arbitrary sisnis ultrafilter, say $\scr W=\scr
  U\sm_n\scr V_n$, and an arbitrary set $A\in\scr W$.  We must find a
  subset of $A$ that is in \scr W and also in \bbb P.  Referring to the
  four clauses in the definition of \bbb P, we see that $A$ already
  satisfies half of the first clause: It has infinitely many infinite
  sections simply because \scr U and the $\scr V_n$'s are
  nonprincipal.  To achieve clause~2 (pairwise disjoint sections), we
  intersect $A$ with a set in \scr W on which $\pi_2$ is one-to-one.
  Such a set exists by part~(2) of Proposition~\ref{sisnis-prop}, and
  the resulting intersection $A'$ satisfies clause~(2) of the
  definition of \bbb P. 

  Next, we shrink $A'$ to an $A''\in\scr W$ satisfying the fourth
  clause in the definition of \bbb P, namely that the $x$-coordinates
  are distinct from the $y$-coordinates of elements of $A''$.  This
  was already done during the proof of \cite[Theorem~2.17]{nbt}, but
  for convenience we repeat the brief argument here.  It suffices to
  find some $B\in\scr U$ such that, for all $n$, we have $B\notin\scr
  V_n$, for then we can take $A''=A'\cap(B\times(\omega-B))$.  Since
  \scr U is distinct from the $\scr V_n$'s, we have, for each $n$,
  some $B_n\in\scr U$ such that $B_n\notin\scr V_n$.  Because \scr U
  is selective and therefore a P-point, it contains a set $B$ such
  that $B-B_n$ is finite for every $n$.  Then, for each $n$, we have
  that $\scr V_n$ contains $\omega-B_n$ and therefore also contains
  its almost-superset $\omega-B$, as required.

Finally, to achieve the remaining half of the first clause (no finite
sections) and the third clause (no elements below the diagonal), we
need only remove finitely many elements from some sections of $A''$;
that will not affect whether those sections are in the ultrafilters
$\scr V_n$ and so the resulting set $A'''$ will be in \scr W.  This
completes the proof that $A$ has a subset $A'''\in\scr W\cap\bbb P$.
\end{pf}

\section{Weak Ramsey Properties}

In this section, we begin the analysis of the connections between the
$(n,T(n))$-weak Ramsey properties in part~(5) of
Propositions~\ref{nbt-prop} and \ref{sisnis-prop} as well as the
three-functions property expressed by part~(3) in these propositions.
We shall study these properties in the context of non-P-points.  This
context makes the weak Ramsey properties quite strong, in the sense
that the next stronger such properties, $(n,T(n)-1)$-weak Ramseyness,
are impossible for non-P-points.

Throughout this section, we assume that we are dealing with a
non-P-point \scr W.  Replacing \scr W by its isomorphic image under a
suitable function to $\omega^2$, we assume further that \scr W is an
ultrafilter on $\omega^2$ such that the first projection
$\pi_1:\omega^2\to\omega$ is neither finite-to-one nor constant on any
set in \scr W.  We refer to such a \scr W as a \emph{non-P-point in
  standard position}.

Notice that, for such a \scr W, every set $A\in\scr W$ has infinitely
many infinite sections, and therefore has a subset in \bbb P.  (The
subset might not be in \scr W.)

We begin by showing that the weak homogeneity in the definition of
$(n,T(n))$-weak Ramseyness necessarily arises from full homogeneity
for each $n$-type.

\begin{notat}
  For any set $S\subseteq\omega^2$ and any $n$-type $\tau$, denote by
  $[S]_\tau$ the set of all those $n$-element subsets of $S$ that
  realize the type $\tau$.
\end{notat}

\begin{prop}            \label{type-homog} 
  Let \scr W be an $(n,T(n))$-weakly Ramsey non-P-point in standard
  position, let $\tau$ be any $n$-type, and let the set
  $[\omega^2]_\tau$ be partitioned into finitely many pieces.  Then
  there is a set $H\in\scr W$ such that $[H]_\tau$ is included in one
  of the pieces.
\end{prop}

\begin{pf}
  Let \scr W, $\tau$, and a partition $\Pi$ of $[\omega^2]_\tau$ into,
  say, $p$ pieces be as in the hypothesis of the proposition.  Define
  a partition $\Pi'$ of $[\omega^2]^n$ into $p+T(n)-1$ pieces by letting the
  first $p$ pieces be those of $\Pi$ and letting the remaining
  $T(n)-1$ pieces be $[\omega^2]_\sigma$ for the $T(n)-1$ $n$-types
  $\sigma$ other than $\tau$.  As \scr W is $(n,T(n))$-weakly Ramsey,
  let $H\in\scr W$ be such that $[H]^n$ meets only $T(n)$ pieces of
  $\Pi'$.  

Now $H$, being in \scr W, has a subset in \bbb P, and, as we noticed
right after Definition~\ref{realize:df}, such a subset contains
realizers for all $n$-types.  Therefore, $[H]^n$ meets all those
pieces of our partition that have the form $[\omega^2]_\sigma$ for
$\sigma\neq\tau$.  That's $T(n)-1$ pieces, so $[H]^n$ can meet at most
one of the remaining pieces of $\Pi'$, which are the original pieces
of $\Pi$.  Therefore, $[H]_\tau$ is included in that single piece of
$\tau$.  
\end{pf}

We shall refer to the conclusion of this proposition as
$\tau$-\emph{homogeneity} or, when we want to refer to all types
$\tau$ together, as \emph{$n$-type homogeneity} for \scr W and for
$H$.  We note that the converse of the proposition is easy when \scr W
has a basis of sets from \bbb P.  For such \scr W, $n$-type
homogeneity implies $(n,T(n))$-weak Ramseyness.  To verify this,
consider any partition of $[\omega^2]^n$ into finitely many pieces and
find, for each $n$-type $\tau$, a homogeneous set $H_\tau\in\scr W$
for that type.  Then all subsets of $\bigcap_\tau H_\tau$ that realize
$n$-types lie in at most $T(n)$ pieces of the original partition,
namely the pieces that contain the sets $[H_\tau]_\tau$.  Finally,
shrink $\bigcap_\tau H_\tau$ to a set in $\scr W\cap\bbb P$, so that
all its $n$-element subsets realize $n$-types.

Notice that the assumption, in the preceding paragraph, that \scr W
has a basis consisting of sets in \bbb P, was used only at the end of
the argument, to ensure that all $n$-element subsets of the
homogeneous set realize some $n$-types.  In some situations, the
property that all $n$-element subsets realize $n$-types can be
obtained from other hypotheses.  The following lemma is a useful
instance of this. It provides, for weakly Ramsey ultrafilters, some
information that would be automatic for \bbb P-generic ultrafilters
and for sisnis ultrafilters, because these are generated by sets in
\bbb P.

\begin{la}              \label{types:la}
If \scr W is an $(n,T(n))$-weakly Ramsey non-P-point in standard
position, then there is a set $P\in\scr W$ such that every $n$-element
subset of $P$ realizes an $n$-type.
\end{la}

\begin{pf}
Partition $[\omega^2]^n$ into $T(n)+1$ pieces by making each of the
$T(n)$ sets $[\omega^2]_\tau$, for $n$-types $\tau$ a piece, and then
adding one more piece containing all the $n$-element sets that don't
realize a type (because they have two elements with the same
$y$-coordinate, or an element whose $x$-coordinate equals another
element's $y$-coordinate, or because an element isn't above the
diagonal).  By hypothesis, there is a set $P\in\scr W$ that meets only
$T(n)$ of these pieces.  Recall that every set in \scr W has a subset
in \bbb P and every set in \bbb P has subsets realizing all types.  So
$[P]^n$ meets all of the pieces of the form $[\omega^2]_\tau$ in our
partition and must therefore miss the one remaining piece, the piece
consisting of $n$-element sets that don't realize types.
\end{pf}

\begin{cor}             \label{down}
If a non-P-point is $(n+1,T(n+1))$-weakly Ramsey, then it is also
$(n,T(n))$-weakly Ramsey.
\end{cor}

\begin{pf}
  We assume, without loss of generality, that \scr W is in standard
  position.  By Lemma~\ref{types:la}, there is a set $A\in\scr W$ all
  of whose $(n+1)$-element subsets realize $(n+1)$-types.  It follows
  immediately that every $n$-element subset of $A$ realizes an
  $n$-type.  This observation allows us to apply the comments
  immediately preceding Lemma~\ref{types:la} without needing the
  assumption that \scr W has a basis of sets in \bbb P.  That
  assumption was needed only to ensure that the final homogeneous set
  can be shrunk so that all its $n$-element subsets realize $n$-types.

Thanks to those comments, it suffices to prove $\tau$ homogeneity for
each $n$-type $\tau$.  Enlarge $\tau$ to an $(n+1)$-type $\tau'$ by
appending $x_{n+1}<y_{n+1}$ after all of the $x$'s and $y$'s
pre-ordered by $\tau$. Notice that, if an $n+1$-element set realizes
$\tau'$ then its first $n$ elements (first in the usual ordering by
$y$-coordinates) form an $n$-element set realizing $\tau$.  Given a
partition $\Pi$ of $[\omega^2]_{\tau}$, form a new partition $\Pi'$ of
$[\omega^2]_{\tau'}$ by putting two sets realizing $\tau'$ into the
same piece of the new partition if their initial $n$-element subsets
(realizing $\tau$) are in the same piece of $\Pi$. By hypothesis, \scr
W contains a set homogeneous for $\Pi'$, and it immediately follows
that this set is also homogeneous for $\Pi$.  (This uses the trivial
fact that every $n$-element set realizing $\tau$ is the initial
$n$-element subset of some $(n+1)$-element set realizing $\tau'$; just
adjoin an $(n+1)\th$ element far beyond the given elements.)
\end{pf}

Our next goal is to establish a connection between weak Ramseyness and
the three-functions property in part~(3) of
Propositions~\ref{nbt-prop} and \ref{sisnis-prop}.  In connection with
the $(2,4)$ in the following proposition, recall that $T(2)=4$.

\begin{prop}            \label{wsel-prop} 
  If a non-P-point \scr W in standard position is $(2,4)$-weakly
  Ramsey, then it has the three-functions property, i.e., every
  function on $\omega^2$ is, on some $H\in\scr W$, either constant or
  one-to-one or $\pi_1$ followed by a one-to-one function.
\end{prop}

\begin{pf}
  Let \scr W be as in the hypothesis of the proposition.  By
  Lemma~\ref{types:la}, let $P\in\scr W$ be a set all of whose
  two-element subsets realize 2-types.

  As a preliminary step, we show that $\pi_1(\scr W)$ is selective.
  Let $g$ be any function on $\omega$; we show that is it one-to-one
  or constant on some set in $\pi_1(\scr W)$.  Let $\tau$ be the type
  given in list form by $x_1<x_2<y_1<y_2$.  Partition the set of pairs
  $\{\sq{a_1,b_1},\sq{a_2,b_2}\}$ into two pieces according to
  whether $g(a_1)=g(a_2)$.  Let $H\in\scr W$ be $\tau$-homogeneous for
  this partition.  Because \scr W is in standard position, $\pi_1$ is
  not finite-to-one on any set in \scr W, so we can arrange that all
  nonempty sections of $H$ are infinite.  It easily follows that $g$
  is constant or one-to-one on $\pi_1(H)\in\pi_1(\scr W)$.

  Now to prove the three-functions property, let $f$ be any function
  on $\omega^2$.  Partition $[\omega^2]^2$ into two pieces, the first
  being $\{\{x,y\}:f(x)=f(y)\}$ and the second being
  $\{\{x,y\}:f(x)\neq f(y)\}$.

Apply $\tau$-type homogeneity for all four 2-types $\tau$ and
intersect the resulting homogeneous sets with each other and with $P$.
The result is a set $H\in\scr W$ such that all its 2-element subsets
realize 2-types and, for each 2-type $\tau$, either $f(a)=f(b)$
whenever $\{a,b\}\subseteq H$ realizes $\tau$ or $f(a)\neq f(b)$
whenever $\{a,b\}$ realizes $\tau$.

Consider first the type given in list form by $x_1=x_2<y_1<y_2$, i.e.,
the type of pairs that lie in a vertical column.  Suppose all pairs
$\{a,b\}\subseteq H$ realizing this type have $f(a)=f(b)$, so $f$ is
constant on vertical columns of $H$.  Then, on $H$, we have
$f=g\circ\pi_1$ for some function $g$ on $\omega$.  Since $\pi_1(\scr
W)$ is selective, it contains a set $K$ on which $g$ is constant or
one-to-one.  Then on $H\cap\pi_1^{-1}(K)\in\scr W$, the function $f$
is either constant (if $g$ is constant on $K$) or $\pi_1$ followed by
a one-to-one function (namely $g$).  So the proposition is true in
this case.

There remains the case that all pairs $\{a,b\}$ realizing the type
$x_1=x_2<y_1<y_2$ have $f(a)\neq f(b)$, so $f$ is one-to-one on
vertical columns of $H$.  We want to show that $f$ is globally
one-to-one on $H$.  We shall do this by considering, one at a time,
the three remaining 2-types, the types realized by pairs not on a
vertical column, and show, in each of the three cases, that $f$ must
take different values at the two points of such a pair.  In each case,
we shall show this by contradiction, assuming that $f$ takes the same
value at the two points of each pair realizing the type, and deducing
that $f$ also takes the same value at two points in a vertical column.

Consider first the 2-type $\tau$ given by $x_1<x_2<y_1<y_2$ and
suppose, toweard a contradiction, that $f(a)=f(b)$ for some, and
therefore for all pairs $\{a,b\}\in[H]_\tau$.  Let $\tau'$ be the
3-type $x_1=x_2<x_3<y_1<y_2<y_3$ and consider any 3-element subset
$\{a,b,c\}$ of $H$ realizing $\tau'$.  (Such a set exists because
every set in \scr W contains realizers for all types.) By inspection
of the types, we see that both $\{a,c\}$ and $\{b,c\}$ realize $\tau$.
therefore we have both $f(a)=f(c)$ and $f(b)=f(c)$.  But $a$ and $b$
are in the same vertical column, so $f(a)\neq f(b)$, and we have a
contradiction.

The remaining two cases are very similar to the one just considered.
When $\tau$ is $x_2<x_1<y_1<y_2$, we take $\tau'$ to be
$x_3<x_1=x_2<y_1<y_2<y_3$, and when $\tau$ is $x_1<y_1<x_2<y_2$, we
take $\tau'$ to be $x_1=x_2<y_1<y_2<x_3<y_3$.  The rest of the
argument is verbatim as in the preceding paragraph.  (The general
recipe for producing $\tau'$ from $\tau$ is to first change the
subscript 2 in $\tau$ to 3 and then insert $x_2$ equal to $x_1$ and
$y_2$ immediately after $y_1$.) 
\end{pf}

The results proved so far in this section give the chain of
implications, in which ``w.r.'' abbreviates ``weakly Ramsey'' and
``3f'' abbreviates ``the three-functions property''.
\[
\dots\to(n+1,T(n+1))\text{-w.r.}\to (n,T(n))\text{-w.r.}\to\dots\to
(2,4)\text{-w.r.} \to3\text{f}.
\]
The question naturally arises whether these implications can be
reversed.  There is a reasonable hope for reversals, by analogy with
what happens for selective ultrafilters.  There, we have easy
implications from the partition properties for larger exponents $n$ to
the partition properties for smaller $n$ and from the partition
property for exponent 2 to selectivity, just as in the chain above.
Kunen's theorem gives reversals for the whole chain by showing that
selectivity implies all the finite-exponent partition properties.  Is
there an analog of Kunen's theorem in the present situation, i.e.,
does the three-functions property (for a non-P-point) imply
$(n,T(n))$-weak Ramseyness for all $n$?  We conclude this section by
showing that the answer is ``no''; in the next section, though, we
will show a way to correct the problem.

\begin{thm}             \label{noreverse} 
  Assume the continuum hypothesis.  There is a non-P-point in standard
  position that satisfies the three-functions property but is not
  $(2,4)$-weakly Ramsey.
\end{thm}

\begin{pf}
   Because the continuum hypothesis is assumed, there is an
  enumeration, in an $\omega_1$-sequence, of all the functions
  $\omega^2\to\omega$; fix such an enumeration
  \sq{f_\alpha:\alpha<\omega_1}.

Using the continuum hypothesis again, fix a selective ultrafilter \scr
U on $\omega$.

Also, fix a symmetric, irreflexive, binary relation $E$ on $\omega$
that makes $\omega$ a copy of the random graph.  This means that,
given any finite list of elements $a_0,\dots,a_{n-1}\in\omega$ and any
subset $S\subseteq n$, there is some $b\in\omega$ distinct from all
the $a_i$'s such that, for all $i<n$, we have $a_iEb$ if and only if
$i\in S$.  In this situation, we say that $b$ realizes the
configuration\footnote{The word ``type'' is often used instead of
  ``configuration''; we prefer the latter here, to avoid any confusion
  with the ``$n$-type'' terminology.}  \sq{a_0,\dots,a_{n-1};S}; we
refer to the $a_i$'s as the parameters of this configuration.  It is
easy to define such a relation $E$ by an inductive construction,
realizing all configurations, one at a time.  It is well known and
easy to prove, by a back-and-forth argument, that any two such
relations $E$ yield isomorphic graphs \sq{\omega,E}.  So it makes
sense to speak of a copy of the (rather than a) random graph.  Call a
subset $X$ of $\omega$ \emph{rich} if it has a subset $Y\subseteq X$
on which the restriction of $E$ is a copy of the random graph.

With these preliminary items available, we are ready to construct the
desired ultrafilter \scr W on $\omega^2$.  It will be generated by
$\aleph_1$ sets $S_\alpha$ indexed by the ordinals $\alpha<\omega_1$,
and subject to the following requirements.
\begin{lsnum}
\item \scr U-almost all sections of $S_\alpha$ are rich, i.e.,
  $\{n\in\omega: S_\alpha(n)\text{ is rich}\}\in\scr U$.
\item The $S_\alpha$'s are decreasing modulo \scr U, i.e., if
  $\alpha<\beta$ then $\{n\in\omega:S_\beta(n)\subseteq
  S_\alpha(n)\}\in\scr U$.
\item The restriction of $f_\alpha$ to $S_{\alpha+1}$ either is
  one-to-one or factors through the first projection, $f_\alpha\restr
  S_{\alpha+1}=g\circ\pi_1$ for some $g:\omega\to\omega$.
\end{lsnum}
Before constructing the sets $S_\alpha$, we verify, on the basis of
these three requirements, that the $S_\alpha$'s together with the sets
$\pi_1^{-1}(X)$ for $X\in\scr U$ generate an ultrafilter
as required in the theorem.

First, we verify that the proposed generators have the finite
intersection property, so they generate a filter.  Consider any
finitely many of the proposed generators, say
$S_{\alpha_1},\dots,S_{\alpha_k}$ and
$\pi_1^{-1}(X_1),\dots,\pi_1^{-1}(X_l)$.  Let $\alpha$ be the largest
of the $\alpha_i$'s and notice that, by requirement~(2), there is a
set $Y\in\scr U$ such that $S_{\alpha_i}(n)\supseteq S_\alpha(n)$ for
all $n\in Y$ and all $i$.  Let $Z$  be the interseection of $Y$ and
all of the $X_j$'s, so $Z\in\scr U$.  Then, whenever $n\in Z$, we have
that the intersection of all of  $S_{\alpha_1},\dots,S_{\alpha_k}$ and
$\pi_1^{-1}(X_1),\dots,\pi_1^{-1}(X_l)$ has the same section at $n$ as
$S_\alpha$ does.  In particular, this section is rich, by
requirement~(1), and, a fortiori, infinite.  This completes the
verifiction that the sets $S_\alpha$ for $\alpha<\omega_1$ and
$\pi_1^{-1}(X)$ for $X\in\scr U$ generate a filter \scr W on
$\omega^2$. 

The preceding argument shows something more, which we record here as a
lemma for future reference.

\begin{la}              \label{Urich} 
  For every set $A\in\scr W$, there is a set $Z\in\scr U$ such that
  $A(n)$ is rich, and therefore infinite, for all $n\in Z$.
\end{la}

Next, we verify that the filter \scr W is an ultrafilter.  Let any
subset $A$ of $\omega^2$ be given, and let $f:\omega^2\to\omega$ be
one-to-one on $A$ and constant on $\omega^2-A$.  This $f$ is
$f_\alpha$ for some $\alpha<\omega_1$, and so it is either one-to-one
or fiberwise constant on $S_{\alpha+1}\in\scr W$ by
requirement~(3). This means that $S_{\alpha+1}$ is (up to perhaps one
point) included in $A$ or in $\omega^2-A$.  Since $A$ was arbitrary,
this proves that \scr W is an ultrafilter.

It now follows immediately from Lemma~\ref{Urich} that \scr W is a
non-P-point in standard position.  

To verify the three-functions property, consider any function
$\omega^2\to\omega$.  It is $f_\alpha$ for some $\alpha<\omega_1$, so
its restriction to $S_{\alpha+1}$ either is constant or factors
through $\pi_1$, by requirement~(3).  Since $S_{\alpha+1}\in\scr W$,
we have the desired conclusion if $f_\alpha$ is constant on this set.
So assume $f_\alpha=g\circ\pi_1$ on $S_{\alpha+1}$.  Because \scr U is
selective, it contains a set $H$ on which $g$ is constant or
one-to-one.  Then $S_{\alpha+1}\cap\pi_1^{-1}(H)$ is a set in \scr W
on which $f_\alpha$ is constant or of the form $g\circ\pi_1$ with $g$
one-to-one.  This completes the verification of the three-functions
property for \scr W.

To show that \scr W is not (2,4)-weakly Ramsey, it suffices, by
Proposition~\ref{type-homog}, to show that \scr W lacks $\tau$-type
homogeneity for some 2-type $\tau$. We use the type given by
$x_1=x_2<y_1<y_2$, i.e., the type of vertical pairs.  These pairs are
partitioned by the edge relation $E$; more precisely, we use the
partition of $[\omega^2]_\tau$ into two pieces, one of which is
$\{\{\sq{a,b},\sq{a,c}\}\in[\omega^2]_\tau:bEc\}$.  For any set
$H\in\scr W$, Lemma~\ref{Urich} implies that $[H]_\tau$ meets both
pieces of this partition.  This completes the verification, on the
basis of requirements~(1)--(3), that \scr W is as required in the
theorem.  

It remains to produce sets $S_\alpha$ satisfying the three
requirements, and we shall do this by recursion on $\alpha$.  We begin
the construction by setting $S_0=\omega^2$; this clearly satisfies
requirement~(1), and the other two requirements are vacuous at this
stage.

We consider next the limit stages of the construction.  Let $\lambda$
be a countable limit ordinal, and suppose, as an induction hypothesis,
that $S_\alpha$ has been defined for all $\alpha<\lambda$, satisfying
requirements~(1)--(3).  Fix an increasing $\omega$-sequence
$0=\alpha(0)<\alpha(1)<\dots$ with supremum $\lambda$.  We shall
define $S_\lambda$ in such a way that \scr U-almost all of its
sections are rich (so that requirement~(1) continues to be satisfied)
and so that $\{n\in\omega:S_\lambda(n)\subseteq
S_{\alpha(i)}(n)\}\in\scr U$ for each $i\in\omega$.  This will ensure
that requirement~(2) holds when $\beta=\lambda$ and $\alpha$ is one of
the $\alpha(i)$'s, but it immediately implies the same for all
$\alpha<\lambda$.  Indeed, for any such $\alpha$, there is an $i$ (in
fact infinitely many of them) with $\alpha<\alpha(i)$, and then,
thanks to the induction hypothesis, we have that
$S_\lambda(n)\subseteq S_{\alpha(i)}(n)\subseteq S_\alpha(n)$ for \scr
U-almost all $n$.  Note that requirement~(3) is exclusively about
successor $S$'s and thus imposes no condition on $S_\lambda$.

As a preliminary normalization, we modify the sequence of sets
\sq{S_{\alpha(i)}:i<\omega}, which is, by requirement~(2) in the
induction hypothesis, decreasing modulo \scr U, to make it literally
decreasing.  That is, we set
\[
S'_i=\bigcap_{j=0}^iS_{\alpha(j)},
\]
and we note that this has not ruined requirement~(1) of the induction
hypothesis; \scr U-almost all sections of $S'_i $ are rich.  Indeed,
for \scr U-almost all $n$, we have that $S_{\alpha(i)}(n)\subseteq
S_{\alpha(j)}(n)$ for all of the (finitely many) $j<i$ and so
$S'_i(n)=S_{\alpha(i)}(n)$. 

We also observe that, since $\alpha(0)=0$, we have
$S'_0=S_{\alpha(0)}=S_0=\omega^2$.  In particular, all sections of
$S'_0$ are rich.

For each $n\in\omega$, let $h(n)$ be the largest $i\leq n$ such that
$S'_i(n)$ is rich. This exists because of the observation in the
preceding paragraph.  Then define
\[
S_\lambda=\{\sq{n,y}\in\omega^2:y\in S'_{h(n)}(n)\}.
\]
The definition of $h$ ensures that $S_\lambda(n)$ is rich for all
$n\in\omega$, so requirement~(1) is satisfied at $\lambda$.
Requirement~(3) says nothing about limit stages, so it remains only to
verify that, for all $i$, \scr U-almost all $n$ satisfy
$S_\lambda(n)\supseteq S_{\alpha(i)}(n)$ for all $i$.  

We remarked, immediately after defining $S'_i$, that it agrees with
$S_{\alpha(i)}$ in \scr U-almost all vertical sections.  So it
suffices to verify that, for each $i$, \scr U-almost all $n$ satisfy
$S_\lambda\subseteq S'_i$.

So fix some $i\in\omega$.  As we pointed out above, \scr U-almost all
$n$ have $S'_i(n)$ rich, and of course, as \scr U is nonprincipal,
\scr U-almost all $n$ have $i\leq n$.  Therefore, \scr U-almost all
$n$ have $h(n)\geq i$.  For these $n$, we have, since the $S'$
sequence is literally decreasing,
\[
S_\lambda(n)=S'_{h(n}(n)\subseteq S'_i(n),
\]
as desired.  This completes the limit stage of our induction.

Finally, we deal with the successor stage.  Suppose, as an induction
hypothesis, that $S_\beta$ has been defined for all $\beta\leq\alpha$
and satisfies requirements~(1)--(3) in this range; our objective is to
define $S_{\alpha+1}$ so as to satisfy all the new instances of
requirements~(1)--(3).  In the case of requirement~(2), it suffices to
ensure that $S_{\alpha+1}(n)\subseteq S_\alpha(n)$ for \scr U-almost
all $n$; the rest of requirement~(2), comparing $S_\beta$ with
$S_{\alpha+1}$ for smaller $\beta$, is then an immediate consequence
via the induction hypothesis comparing $S_\beta$ with $S_\alpha$.

To avoid unnecessary clutter, let us take advantage of the fact that,
in this step of the construction $\alpha$ is fixed, so we can write
simply $f$ for $f_\alpha$. We write $f^n$ for the restriction of $f$
to the $n\th$ column in $\omega^2$, regarded as a function on
$\omega$; that is, $f^n(y)=f(\sq{n,y})$.

 Our goal is thus to produce a set
$S=S_{\alpha+1}$ with the following properties, corresponding to the
three requirments above.
\begin{lsnum}
  \item \scr U-almost all sections $S(n)$ are rich.
\item $S(n)\subseteq S_\alpha(n)$ for \scr U-almost all $n$.
\item The restriction of $f$ to $S$ either is one-to-one or factors
  through $\pi_1$.
\end{lsnum}

Call a number $n$ good (for this argument) if $f^n$ is constant on
some rich subset of $S_\alpha(n)$; otherwise, call $n$ bad.  If \scr
U-almost all $n$ are good, then our task is easy.  For each good $n$,
choose some rich $R_n\subseteq S_\alpha(n)$ on which $f^n$ is
constant, and set 
\[
S=\{\sq{n,y}\in\omega^2:n\text{ is good and }y\in R_n\}.
\]
This has rich sections at all good $n$, so we have (1), while (2)
holds because $R_n\subseteq S_\alpha(n)$ and (3) holds because $f^n$
is constant on $R_n$.

It remains to treat the case that \scr U-almost all $n$ are bad.  In
this case, we shall construct $S$ one element at a time, making sure
that they all map to different values under $f$, so that $f$ will be
one-to-one on $S$.  The elements that we put into $S$ will all be of
the form \sq{n,y} with $n$ bad and $y\in S_\alpha(n)$.  This will
ensure that (2) holds.  The sections of $S$ at good $n$ will be empty,
but this does no harm to (1) because \scr U-almost no $n$ are good.
We shall also ensure that every section of $S$ at a bad $n$ is not
only rich but a copy of the random graph.  (Recall that a rich set is
a superset of a copy of the random graph.) We may also assume that the
nonempty sections of $S_\alpha$ are copies of the random graph; simply
shrink each rich section to a copy of the random graph, and remove all
non-rich sections.

To inductively produce $S$, we first make a list of conditions that
should be satisfied by our construction.  These conditions are
represented formally by tuples of the form
\sq{k_0,k_1,\dots,k_{l-1};Q} where $k_0<k_1<\dots<k_{l-1}$ are natural
numbers and $Q\subseteq l$.  The meaning of this tuple is: ``If the
elements put into $S$ at the $k_0\th,k_1\th,\dots,k_{l-1}\th$ steps of
the construction exist\footnote{Some steps won't put any elements into
  $S$, so ``exist'' is not a vacuous requirement here.} and are in the
same column, say they are $\sq{n,y_1},\sq{n,y_2},\dots,\sq{n,y_{l-1}}$
for some bad $n$, then put into $S$ an element \sq{n,z} such that, for
all $i<l$, we have $zEy_i$ if and only if $i\in Q$.''  Thus, the tuple
\sq{k_0,k_1,\dots,k_{l-1};Q} requests fulfillment of one instance of
the condition that $S(n)$ be a copy of the random graph.  There are
only countably many such tuples, so we can enumerate them as an
$\omega$-sequence.  Fix such an enumeration in which, for all
$q\in\omega$, the $q\th$ tuple has all of its $k_i$'s smaller than $q$
(if necessary, repeat the vacuous tuple, where $l=0$, numerous times).

We now explain the $q\th$ step of the construction of $S$.  Consider
the $q\th$ condition in our enumeration, say
\sq{k_0,k_1,\dots,k_{l-1};Q}, and suppose that, as in the meaning of
this condition explained above, the $k_i$'th step (which has already
been done, because of the way we arranged the enumeration of
conditions) put \sq{n,y_i} into $S$ for each $i<l$.  We wish to adjoin
some \sq{n,z} to $S$ subject to two desiderata.  First, it should do
what the condition \sq{k_0,k_1,\dots,k_{l-1};Q} requests; $zEy_i$
should hold when $i\in Q$ and fail when $i\in l-Q$.  Second,
$f(\sq{n,z})$ should be different from $f(a)$ for all of the finitely
many elements $a$ already put into $S$ during previous steps.  Let us
call these finitely many values $f(a)$ the forbidden values.  

From now on, we work within $S_\alpha(n)$, which we recall is, with
the binary relation $E$, a copy of the random graph.  We seek an
element $z$ that has the correct configuration relative to the
$y_i$'s, as specified by \sq{k_0,k_1,\dots,k_{l-1};Q}, and such that
$f^n(z)$ is different from the forbidden values.  

For each forbidden value $v$, consider the set $B_v$ obtained by
removing from $(f^n)^{-1}(\{v\})$ any $y_i$'s that happen to lie in
$(f^n)^{-1}(\{v\})$.  Since $f^n$ is constant on $B_v$ (with value
$v$) and since $n$ is bad, we know that $B$ is not a copy of the
random graph.  So we can fix a finite subset $F_v$ of $B_v$ and a
configuration $C_v$ with respect to $F_v$ that is not realized by any
element of $B_v$.  Combine all these configurations $C_v$, and also
the configuration that \sq{k_0,k_1,\dots,k_{l-1};Q} tells us to
realize, into a single large but finite configuration $C$, relative to
all the members of the $F_v$'s and all the $y_i$'s.  Since $S_\alpha$
is a copy of the random graph, it contains an element $z$, distinct
from all members of the $F_v$'s and from the $y_i$'s, and realizing
$C$.  Then, for each forbidden value $v$, we have that $z$ realizes
$C_v$ and is therefore not in $B_v$.  It is also not among the $y_i$'s
that were removed from $(f^n)^{-1}(\{v\})$ when we defined $B_v$, so
$z\notin(f^n)^{-1}(\{v\})$.  That is, $f^n(z)$ is not the forbidden
value $v$.  Furthermore, $z$ realizes the configuration requested by
\sq{k_0,k_1,\dots,k_{l-1};Q}.  Therefore, $z$ satisfies all our
desiderata, and we adjoin \sq{n,z} as the next element of $S$.  This
completes the $q\th$ step of the construction of $S=S_{\alpha+1}$.  So
we have obtained $S_{\alpha+1}$ with the required properties, and the
proof of the theorem is complete.
\end{pf}

The theorem just proved shows that we do not have a perfect analog of
Kunen's results for selective ultrafilters.  The three-functions
property, which is the analog of selectivity (a two-functions
property) in our non-P-point situation, does not imply the strongest
possible Ramsey properties.  There remain at least three natural
questions. 
\begin{ls}
  \item Does the three-functions property imply weaker Ramsey
properties, say $(n,t)$-weak Ramseyness for some $t>T(n)$? 
\item  Do any of the $(n,T(n))$-weak Ramsey properties imply
other such properties with larger $n$?
\item Can the three-functions property be combined with some other
  (reasonable) property to imply weak Ramsey properties?
\end{ls}

The second and third of these questions will be answered in the next
section.  As for the first, the following remark sketches a negative
answer.

\begin{rmk} \label{cbly-random} The random graph can be viewed as a
  random edge 2-coloring of the complete graph on $\aleph_0$ vertices.
  Edges of the random graph are colored red, and edges of the complete
  graph that are not in the random graph are colored green.  There is
  an entirely analogous random edge $\aleph_0$-coloring of the
  complete graph on $\aleph_0$ vertices. It defining property is that,
  given any finite set $F$ of vertices and any function $f$ assigning
  to each vertex $v\in F$ one of the colors, there is a vertex
  $z\notin F$ whose edge to any $v\in F$ has the color $f(v)$.

The proof of Theorem~\ref{noreverse} can be carried out essentially
unchanged but with this random $\aleph_0$-coloring in place of the
random graph.  The result is a non-P-point \scr W in standard
position, enjoying the three-functions property, but with the
following strong negative partition property.  For the 2-type $\tau$
given by $x_1=x_2<y_1<y_2$, there is a partition of $[\omega^2]_\tau$
into infinitely many pieces such that, for every $H\in\scr W$, all of
the pieces meet $[H]_\tau$.  

In particular, by merging some of the pieces of this partition, we can
get, for any finite number $t$, a partition of $[\omega^2]_\tau$ into
$t$ pieces such that all pieces meet $[H]_\tau$ for all $H\in\scr W$.
Combining this partition of $[\omega^2]_\tau$ with the three other
pieces $[\omega^2]_\sigma$ for 2-types $\sigma\neq\tau$, we get a
counterexample showing that \scr W is not $(2,t+2)$-weakly Ramsey. 
\end{rmk}

\section{Conservativity}

In this section, we recall the notion of conservative elementary
extensions, introduced in the context of models of arithmetic by
Phillips \cite{phillips}; we explain its connection with weak Ramsey
properties; and we show that it, when combined with the
three-functions property, implies $(n,T(n))$-weak Ramseyness for all
$n\in\omega$. 

To avoid excessive repetition, we refer the reader to \cite{model-th}
for some of the results that we shall need and we give only a short
summary here.

We adopt the convention that, for a structure \ger A, its underlying
set (also called its domain or its universe) is denoted by $|\ger A|$.

\begin{df}              \label{cons}
Let \ger A be a structure for a first-order language, and let \ger B
be an elementary extension of \ger A.  Then \ger B is a
\emph{conservative extension} of \ger A if, whenever $X$ is a
parametrically definable subset of \ger B, then $X\cap|\ger A|$ is a
parametrically definable subset of \ger A.
\end{df}

This concept makes good sense in the context of general model theory.
In fact, it can be used to characterize stable theories as those
theories $T$ such that all elementary extensions of models of $T$ are
conservative extensions.  We shall, however, use only the special case
where the models are elementary extensions of the standard model \ger
N of full arithmetic.  By full arithmetic, we mean the language that
has function and relation symbols for all of the functions and
relations (of arbitrary finite arity) on the set of natural numbers.
\ger N is the model with underlying set $\omega$ and with all the
symbols having the obvious meanings.  In fact, we shall be concerned
only with ultrapowers $\scr U\pr\ger N$ of the standard model \ger N.  

Since all subsets of \ger N are definable in the language of full
arithmetic, all elementary extensions of \ger N are conservative
extensions. 

If \scr U is an ultrafilter on $\omega$, then $\scr U\prn$ is
generated by a single element, namely the equivalence class modulo
\scr U of the identity function, $[\text{id}]_{\scr U}$.  Indeed,
every element $[f]_{\scr U}$ of $\scr U\prn$ is
${}^*\!f([\text{id}]_{\scr U})$, where, as is customary in nonstandard
analysis, ${}^*\!f$ is the function on $\scr U\prn$ denoted there by
the function symbol that denotes $f$ in \ger N.

If \scr U is an ultrafilter on a countably infinite set $S$ other than
$\omega$, it is still the case that $\scr U\prn$ is generated by a
single element; the equivalence class modulo \scr U of any bijection
between $S$ and $\omega$ will do.

If \scr U is an ultrafilter on $S$ and $f$ is a function with domain
$S$, then $f$ induces an elementary embedding $f_*:f(\scr
U)\prn\to\scr U\prn$, namely the function that sends each $[g]_{f(\scr
  U)}$ to $[g\circ f]_{\scr U}$.  We sometimes identify $f(\scr
U)\prn$ with its image under this embedding.

\begin{df}
If $\scr U\prn$ is a conservative extension of the image under $f_*$
of $f(\scr U)\prn$, then we call $f$ a \emph{conservative map} on \scr
U. 
\end{df}

We shall be particularly interested, for reasons to be explained later,
in the situation where \scr U is an ultrafilter on $\omega^2$ and $f$
is the projection $\pi_1$ to the first cooordinate.

The key property of conservative extensions for our purposes is the
following result, which is part of Theorem~3 in \cite{model-th}.

\begin{prop} \label{unique} Suppose that \ger B and \ger C are
  elementary extensions of a model \ger A of full arithmetic and that
  \ger B is a conservative extension of \ger A.  Then there is, up to
  isomorphism, only one amalgamation of \ger B and \ger C over \ger A
  in which all the elements of $|\ger B|-|\ger A|$ are above (with
  respect to the nonstandard extension ${}^*\!{<}$ of the standard order
  on $\omega$) all elements of $|\ger C|$.
\end{prop}

The existence of such amalgamations is established by a fairly easy
compactness argument; see Theorem~2(b) in \cite{model-th}.  The
important part of Proposition~\ref{unique} is the uniqueness.

The relevance of amalgamations for our purposes is the following
connection with weak Ramsey properties; it is Theorem~5 of
\cite{model-th}.

\begin{prop}            \label{amal-wkR}
  An ultrafilter \scr U on $\omega$ is $(n,t)$-weakly Ramsey if and
  only if there are, up to isomorphism, at most $t$ ways to amalgamate
  $n$ copies of $\scr U\prn$ with a specified ordering of the $n$
  copies of the generator $[\text{id}]_{\scr U}$.
\end{prop}

These propositions allow us to prove the first main result of this
section.

\begin{thm}             \label{cons2wkR}
Let \scr W be a non-P-point in standard position.  Assume that \scr W
has the three-functions property and that $\pi_1$ is a conservative
map on \scr W.  Then \scr W is $(n,T(n))$-weakly Ramsey for all
$n\in\omega$. 
\end{thm}

\begin{pf}
We write \scr U for the ultrafilter $\pi_1(\scr W)$ on $\omega$.

Thanks to the three-functions property, the ultrapower $\scr W\prn$
has exactly three elementary submodels, namely the standard model \ger
N, the whole model $\scr W\prn$, and the copy of $\scr U\prn$ induced
by the projection $\pi_1$ from \scr W to \scr U.  The copy of $\scr
U\prn$ is generated by $[\pi_1]_{\scr W}$.  The whole model $\scr
W\prn$ is generated by the equivalence class, modulo \scr W, of any
bijection $\omega^2\to\omega$.  It is also generated by $[\pi_2]_{\scr
  W}$ because $\pi_2$ is one-to-one on a set in \scr W and therefore
coincides, on a possibly smaller set in \scr W, with a bijection.  

Temporarily, consider a single $n$-type $\tau$.  The notion of an
$n$-tuple of elements of $\omega^2$, with $y$-coordinates in
increasing order, realizing $\tau$ can, like any relation on $\omega$,
be canonically extended to any elementary extension of \ger N.  

In particular, suppose we have an amalgamation of $n$ copies of $\scr
W\prn$, with their generators, copies of $[\pi_2]$, in a specified
order.  Let the $n$ copies of $\scr W\prn$ be listed in the order of
their copies of $[\pi_2]$ in the amalgamation, and write $[f]^i$ for
the image, in the amalgamation, of an element $[f]$ of the $i\th$ copy
of $\scr W\prn$.  So our convention for the order of listing these
copies ensures that $[\pi_2]^1<[\pi_2]^2<\dots<[\pi_2]^n$ in the
amalgamation.  Then the elements $[\text{id}]^i$, which are pairs in
the amalgamation, realize $\tau$ if and only if the list form of
$\tau$ is satisfied when the $x_i$ are interpreted as $[\pi_1]^i$ and
the $y_i$ as $[\pi_2]^i$.

We now check how many amalgamations there are, of $n$ copies of $\scr
W\prn$, in which the generators $[\text{id}]^i$ realize $\tau$.  We
build such an amalgamation by starting with the standard model \ger N
and extending it in as many steps as there are inequivalent elements
in the pre-order $\tau$; we go through these elements in increasing
order according to $\tau$.  

At a step corresponding to an equivalence class of $x_i$'s in $\tau$,
we must amalgamate the model \ger M already constructed in previous
steps with $\scr U\prn$, identifying only the standard parts, and
putting all nonstandard elements of $\scr U\prn$ above all elements of
\ger M.  Since $\scr U\prn$ is, like any model of full arithmetic, a
conservative extension of \ger N, there is only one way, up to
isomorphism, to perform this amalgamation.

At a step corresponding to a $y_i$ in $\tau$, we must amalgamate the
model \ger M already constructed with $\scr W\prn$, identifying the copy
of $\scr U\prn$ in this $\scr W\prn$ with the copy already amalgamated
into \ger M at the earlier step corresponding to $x_i$, and putting
all elements of $|\scr W\prn|-|\scr U\prn|$ above all elements of \ger
M.  Again, there is only one way to perform this amalgamation, because
$\scr W\prn$ is, by hypothesis, a conservative extension of the
submodel identified with $\scr U\prn$.  

The preceding two paragraphs show that there is only one way, up to
isomorphism,  to amalgamate $n$ copies of $\scr W\prn$ with the
generators in a specified order and realizing $\tau$.  Since this
holds for each $n$-type $\tau$, and since there are $T(n)$ $n$-types,
we conclude that there are only $T(n)$ ways to amalgamate $n$ copies
of $\scr W\prn$ with the generators in a specified order.  By
Proposition~\ref{amal-wkR}, it follows that \scr W is
$(n,T(n))$-weakly Ramsey.
\end{pf}

Theorem~\ref{cons2wkR} shows that conservativity is a sufficient
condition to add to the three-function property and produce all the
weak Ramsey properties enjoyed by sisnis ultrafilters and \bbb
P-generic ultrafilters.  The question naturally arises whether
conservativity is necessary for this purpose.  Theorem~\ref{noreverse}
shows that some additional condition is needed, but might
conservativity be excessive? Might a weaker additional condition
suffice?  The next theorem answers these questions negatively.  It
shows that conservativity is really needed.

\begin{thm}             \label{wR2cons}
  Suppose \scr W is a $(2,4)$-weakly Ramsey non-P-point in standard
  position.  Then $\pi_1$ is a conservative map on \scr W.
\end{thm}

\begin{pf}
As before, we write \scr U for $\pi_1(\scr W)$.  As a first step
toward the proof, we analyze the definition of ``$\pi_1$ is a
conservative map on \scr W'' in order to replace it with an equivalent
condition of a combinatorial, rather than model-theoretic, flavor.  

Notice first that, when considering a parametrically definable subset
$X$ of $\scr W\prn$, we may assume without loss of generality that the
only parameter used in the definition is the pair $[\text{id}]_{\scr
  W}$.  This is because any other parameter $[f]_{\scr W}$ can be
defined from $[\text{id}]_{\scr W}$.  Furthermore, although the
definition could, a priori, look like 
\[
X=\{z\in\scr W\prn:\phi([\text{id}]_{\scr W},z)\} 
\]
for an arbitrary formula $\phi$, we may assume without loss of
generality that it has the form 
\[
X=\{z\in\scr W\prn:{}^*\!R([\pi_1]_{\scr W},[\pi_2]_{\scr W},z)\} 
\]
for some ternary relation $R$ on $\omega$.  This is because we are
working in full arithmetic, so any formula $\phi$ applied to a pair
and a single element is equivalent, in \ger N and therefore in any
elementary extension, to an atomic ternary relation.  

So let us consider an arbitrary $X$ of this form, obtained from some
arbitrary ternary relation $R$.  The intersection of $X$ with the
elementary submodel $\scr U\prn$ is, taking into account the
identification via $(\pi_1)_*$, 
\[
Y=\{[f]_{\scr U}: \scr W\prn\models{}^*\!R([\pi_1]_{\scr
  W},[\pi_2]_{\scr W}, [f\circ\pi_1]_{\scr W})\}.
\]

Our goal is to prove that this $Y$ is parametrically definable in $\scr
U\prn$.  As in the case of $X$, if there is such a definition, there
will be one whose only parameter is $[\text{id}]_{\scr U}$ and indeed
one of the form 
\[
Y=\{[f]_{\scr U}:\scr U\prn\models {}^*\!S([\text{id}]_{\scr
  U},[f]_{\scr U})\}
\] 
for some binary relation $S$ on $\omega$.  

Taking into account the definition of how relation symbols are
interpreted in ultrapowers, we find that what must be proved is the
following.  For every ternary relation $R$ on $\omega$, there exists a
binary relation $S$ on $\omega$ such that, for all functions
$f:\omega\to\omega$, we have
\[
\{\sq{x,y}\in\omega^2:R(x,y,f(x))\}\in\scr W\iff
\{x\in\omega:S(x,f(x))\}\in\scr U. 
\]

Now that we have a combinatorial version of the desired conclusion, we
work toward deducing this version form the assumption that \scr W is a
(2,4)-weakly Ramsey non-P-point in standard position.  In fact, we
will use (2,4)-weak Ramseyness only to obtain $\tau$-homogeneity for
the 2-type $\tau$ given by $x_1=x_2<y_1<y_2$.  

Let an arbitrary ternary relation $R$ be given.  Associate to each
pair $\sq{x,y}\in\omega^2$ the function $C_{\sq{x,y}}:\omega\to2$ that
sends any $z\in\omega$ to 1 if $R(x,y,z)$ and to 0 otherwise.
Partition $[\omega^2]_\tau$ into two pieces, putting
$\{\sq{x,y},\sq{x,y'}\}$ (where $y<y'$ by our usual convention) into
the first piece if $C_{\sq{x,y}}$ lexicographically precedes
$C_{\sq{x,y'}}$ and into the second piece otherwise.  By hypothesis,
there is a set $H\in\scr W$ such that $[H]_\tau$ lies entirely in one
of the two pieces.  

Suppose $[H]_\tau$ is included in the first piece of our
partition. (The alternative possibility, that it is included in the
second piece, is handled by an entirely analogous argument.)  We may
assume all nonempty sections of $H$ are infinite, since removing any
finite sections only deletes a set not in \scr W and thus changes none
of the properties we have for $H$.  For each $x\in\pi_1(H)$, the
sequence of functions \sq{C_{\sq{x,y}}:y\in H(x)} is lexicographically
increasing.  Any such sequence eventually stabilizes componentwise.
That is, for each $z\in\omega$, there is some $N_z$ such that
$C_{\sq{x,y}}(z)$ is independent of $y$ once $y\geq N_z$.  To see
this, argue by induction on $z$.  For $z=0$, the lexicographic
ordering forces the values of $C_{\sq{x,y}}(0)$ to never decrease as
$y$ increases, so they are either all 0, or, once one of them is 1,
all the later ones, for larger $y$, are also 1.  Once the values for
$z=0$ have stabilized, the values for $z=1$ can never decrease, so
these too must stabilize.  And so on; once the values for all $z<k$
have stabilized, the values for $z=k$ can no longer decrease, so they
also stabilize.  

Now define $S$ by putting a pair \sq{x,z} into $S$ if and only if the
eventual, stable value of $C_{\sq{x,y}}(z)$ for all sufficiently large
$y$ is 1. We claim that this $S$ works.  Let an arbitrary
$f:\omega\to\omega$ be given.

Suppose first that $\{x\in\omega:S(x,f(x))\}\in\scr U$, and let $B$
denote this set in \scr U.  For each $x\in B$ we have, by definition
of $S$, that $R(x,y,f(x))$ holds for all sufficiently large $y\in H(x)$.
Thus, the set $\{\sq{x,y}\in\omega^2:R(x,y,f(x))\}$ includes the
intersection of $H$, $\pi_1^{-1}(B)$, and a set of the form
$\{\sq{x,y}:y>M(x)\}$ for some function $M$.  All three of these are
in \scr W, the last because of standard position: $\pi_1$ is not
finite-to-one on any set in \scr W.  Therefore the intersection is in
\scr W, as required.

The remaining case, that $\{x\in\omega:S(x,f(x))\}\notin\scr U$ is
handled the same way, using, in place of $R$ and $S$, their
negations.  
\end{pf}

\begin{cor}             \label{alln}
For non-P-points, the properties of $(n,T(n))$-weak Ramseyness for
different $n\geq2$ are all equivalent.
\end{cor}

\begin{pf}
We already know, from Corollary~\ref{down}, that these weak Ramsey
properties for larger $n$ imply the proiperties for smaller $n$.  For
the converse, we assume that the ultrafilter is in standard position;
this can be arranged by applying a suitable bijection and thus entails
no loss of generality.  Then any of these weak Ramsey properties
implies $(2,4)$-weak Ramseyness (by Corollary~\ref{down}), which in
turn implies both the three-functions property (by
Proposition~\ref{wsel-prop}) and conservativity of $\pi_1$ (by
Theorem~\ref{wR2cons}).  These, in turn, imply $(n,T(n))$-weak
Ramseyness for all $n$ (by Theorem~\ref{cons2wkR}).
\end{pf}

Summarizing, we have, for non-P-points in standard position, the
equivalence of all the $(n,T(n))$-weak Ramsey properties and the
conjunction of the three-functions property with conservativity of
$\pi_1$.  For non-P-points not in standard position, the only change
that is needed is that conservativity applies not to $\pi_1$ but to
any function $p$ that is neither finite-to-one nor constant on any set
of the ultrafilter.  (The three-functions property ensures that $p$ is
essentially unique.)

\section{Infinitary Partition Relations and Complete Combinatorics}

In the preceding sections, we have dealt only with finitary partition
relations.  In the present section, we turn to infinitary partition
relations enjoyed by \bbb P-generic ultrafilters and by sisnis
ultrafilters.  By analogy with Mathias's results for selective
ultrafilters, parts~(2) and (3) of Proposition~\ref{sel-ramsey}, and
thinking of \bbb P in our situation as being the analog of
$[\omega]^\omega$ in Mathias's situation, we might hope that our
ultrafilters \scr W enjoy a partition relation of the following sort:
Whenever \bbb P is partitioned into two nice pieces, then there is
some $A\in\scr W$ all of whose subsets in \bbb P lie in the same
piece.  Here, ``nice'' could mean analytic/coanalytic, or, in the case
of the L\'evy-Mahlo model, it could mean $HOD\bbb R$.

Unfortunately, such a partition relation is extremely false.  It is
possible to partition \bbb P into continuum many Borel pieces, all of
which are dense in the forcing notion  \bbb P.  To see this, we extend
the notion of $n$-types (Definition~\ref{type:df}) to $\omega$-types.

\begin{df}              \label{omega-type} 
  An $\omega$-\emph{type} is a linear pre-order of the infinite set of
  formal symbols $x_1,x_2,\dots$ and $y_1,y_2,\dots$ such that
  \begin{ls}
    \item $y_1<y_2<\dots$,
\item each $x_i$ precedes the corresponding $y_i$,
\item each equivalence class in the pre-order consists of either a
  single $y_i$ or infinitely many $x_i$'s, 
\item there are infinitely many equivalence classes of $x$'s, and
\item the induced linear order of the equivalence classes has
  order-type $\omega$.
  \end{ls}
\end{df}

The intention here is that an $\omega$-type describes the
order-relationships between the $x$ and $y$ coordinates of the points
in an element of \bbb P.  Recall that the definition of \bbb P
requires that, if $A\in P$, then all the points in $A$ have distinct
$y$-coordinates; as before, we adopt the convention of thinking of the
points in $A$ as listed in order of increasing $y$-coordinates.

\begin{df}
  The $\omega$-type realized by an element $A$ of \bbb P is the
  pre-order consisting of exactly those inequalities between the
  formal symbols $x_i$ and $y_j$ that hold when $A$ is listed as
  $\{\sq{a_i,b_i}:i\in\omega\}$ in order of increasing $b_i$'s and
  then each $x_i$ is interpreted as denoting $a_i$ and each $y_j$ is
  interpreted as denoting $b_j$.
\end{df}

The definition of \bbb P easily implies that every $A\in\bbb P$
realizes a (unique) $\omega$-type.  

\begin{notat}
  Let $\bbb P_\tau$ be the subset of \bbb P consisting of those
  elements of \bbb P that realize the $\omega$-type $\tau$.
\end{notat}

Given any $A\in\bbb P$ and any $\omega$-type $\tau$, it is easy to
construct a subset $B\subseteq A$ in \bbb P (i.e., an extension of
$A$ in the forcing notion \bbb P) that realizes $\tau$.  That is, each
$\bbb P_\tau$ is dense in \bbb P.  It is easy to check also that each
$\bbb P_\tau$ is a Borel set.  So we have, as claimed, a partition of
\bbb P into continuum many Borel sets, all of which are dense in \bbb
P.

Although this result constitutes a strong counterexample to natural
partition relations for \bbb P, it also suggests a way around the
problem.  Each $\bbb P_\tau$ is a notion of forcing equivalent to \bbb
P, and we might hope for a partition relation satisfied by one of
these notions of forcing.  Recall Remark~\ref{forcing-equiv}, where we
pointed out that an infinite-exponent partition relation can hold for
a notion of forcing while failing for an equivalent notion.  Perhaps
this happens here.

In fact, the next theorem shows that this happens for every
$\omega$-type. 

\begin{thm}             \label{analytic}
Let $\tau$ be an $\omega$-type, and let $\bbb P_\tau$ be partitioned
into an analytic subset and its complement.  
\begin{lsnum}
  \item There is a set $H\in\bbb P_\tau$ such that all its subsets in
    $\bbb P_\tau$ lie in the same piece of the partition.
\item Any sisnis ultrafilter on $\omega^2$, contains an $H$ such that
  all its subsets in $\bbb P_\tau$ lie in the same piece of the partition.
\item Any \bbb P-generic ultrafilter on $\omega^2$, contains an
  $H\in\bbb P_\tau$ such that all its subsets in $\bbb P_\tau$ lie in
  the same piece of the partition.
\end{lsnum}
\end{thm}

Parts~(1) and (3) were proved for a particular $\omega$-type $\tau$ by
Dobrinen in \cite{dobrinen} using an entirely different method, based
on Todorcevic's theory of topological Ramsey spaces.  It is very
likely that her method can be applied to arbitrary $\omega$-types, not
just the one she used in \cite{dobrinen}.

\begin{pf}
  The main work in this proof is to establish part~(2) of the theorem;
  afterward, parts~(1) and (3) will follow fairly easily.
  Fortunately, the main work was already done in \cite{sel-homog},
  specifically in proving Theorem~7 of that paper.  So our first task
  here is just to show how (2) follows from that theorem.  This
  argument parallels part of the proof of Theorem~2.17 in \cite{nbt},
  which also relied on the same result from \cite{sel-homog}.

We begin by stating, in the next lemma, the result from
\cite{sel-homog}; afterward, we shall show how part~(2) of the present
theorem follows from it.

\begin{la}[Theorem~7 of \cite{sel-homog}] \label{old} Assume that
  selective ultrafilters $\scr D(s)$ have been assigned to all finite
  subsets $s$ of $\omega$, and assume that every two of these
  ultrafilters are either equal or not isomorphic.  Let \scr X be an
  analytic subset of $[\omega]^\omega$.  Then there is a function $Z$
  assigning, to each ultrafilter \scr D that occurs among the $\scr
  D(s)$'s, some element $Z(\scr D)\in\scr D$ such that \scr X
  contiains all or none of the infinite sets
  $\{z_0<z_1<z_2<\dots\}\in[\omega]^\omega$ that satisfy $z_n\in
  Z(\scr D(\{z_0,\dots,z_{n-1}\}))$ for all $n\in\omega$.
\end{la}

We emphasize that, if the same ultrafilter \scr D occurs as $\scr
D(s)$ for several sets $s$, then a single set $Z(\scr D)$ is assigned
to it by $Z$, not a possibly different set for each occurrence.

Using this lemma, we proceed with the proof of part~(2) of our
theorem.  Let $\scr W=\scr U\sm_n\scr V_n$ be a sisnis ultrafilter on
$\omega^2$, so \scr U and all of the $\scr V_n$ are non-isomorphic
selective ultrafilters on $\omega$.  Let $\tau$ be an $\omega$-type,
and let $\bbb P_\tau$ be partitioned into an analytic piece \scr Y and
its complement.  

There is a natural bijection $\phi$ from $[\omega]^\omega$ onto $\bbb
P_\tau$, defined as follows.  Given a set
$\{z_0<z_1<z_2<\dots\}\in[\omega]^\omega$, assign the value $z_i$ to
the formal variables in the $i\th$ equivalence class\footnote{Here and
  below, our enumeration of the equivalence classes begins with 0. So
  the $i\th$ equivalence class is the one with exactly $i$ strict
  predecessors.} of the pre-order $\tau$.  (So $z_i$ becomes the value
of either a single $y_j$ or infinitely many $x_k$'s).  For each
$j\in\omega$, the values assigned to $x_j$ and $y_j$ determine a point
in $\omega^2$, and we let $\phi(\{z_0,z_1,\dots\})$ be the set of
these points.  Because $\tau$ is an $\omega$-type, this set is in
$\bbb P_\tau$.

Let $\scr X=\phi^{-1}(\scr Y)$.  Since $\phi$ is clearly continuous
and since \scr Y is analytic, \scr X is an analytic subset of
$[\omega]^\omega$, as required by the hypothesis of Lemma~\ref{old}.

To apply the lemma, we must still define appropriate selective
ultrafilters $\scr D(s)$; these will be chosen from among the
ultrafilters \scr U and $\scr V_n$ that produced the sisnis
ultrafilter \scr W.  The scheme for associating these ultrafilters to
finite subsets $s$ of $\omega$ is as follows.  Suppose
$s=\{z_0<z_1<\dots<z_{k-1}\}$.  Assign the value $z_i$ to the formal
variables in the $i\th$ equivalence class of the pre-order $\tau$, for
each $i<k$.  (This is exactly like the definition of $\phi$ above,
except that, because $s$ is finite, only $k$ equivalence classes of
variables get values.)  Consider equivalence class number $k$, the
first one not assigned a value here.  If it is an equivalence class of
$x_j$'s, then let $\scr D(s)=\scr U$.  If, on the other hand, it
consists of a single $y_j$, then, since the corresponding $x_j$
precedes $y_j$ in $\tau$, $s$ has assigned a value $v$ to $x_j$; set
$\scr D(s)=\scr V_v$.

Since \scr U and all the $\scr V_n$ are non-isomorphic selective
ultrafilters, we have satisfied the hypotheses of Lemma~\ref{old}, so
we obtain a function $Z$ as described there.  It remains to chase
through all the relevant definitions to see what the conclusion of
Lemma~\ref{old} tells us in this situation.  

That conclusion concerns sets $\{z_0<z_1<\dots\}$ such that each $z_n$
is a member of $Z(\scr D(\{z_0,\dots,z_{n-1}\}))$.  By our choice of \scr D's
this means that, when we compute $\phi$ of such a set, the values
assigned to the $x_j$'s are in $Z(\scr U)$ and the value assigned to
any $y_j$ is in $Z(\scr V_v)$ where $v$ is the value of the
corresponding $x_j$.  This means that $\phi(\{z_0<z_1<\dots\})$ is a
subset of
\[
H=\{\sq{a,b}\in\omega^2:a\in Z(\scr U)\text{ and }b\in Z(\scr V_a)\}.
\]
This $H$ is in $\scr W=\scr U\sm_n\scr V_n$ because each $Z(\scr D)$
is in the corresponding \scr D.  Furthermore, any subset of $H$ of
type $\tau$ is $\phi(\{z_0<z_1<\dots\})$ for some $\{z_0<z_1<\dots\}$
as in the conclusion of Lemma~\ref{old}.  The lemma tells us that
either all or none of these sets  $\{z_0<z_1<\dots\}$ are in \scr X,
and, in view of our choice of \scr X, this means that all or none of
the subsets of $H$ of type $\tau$ are in \scr Y.  This completes the
proof of part~(2) of the theorem.

We turn next to part~(1).  We first prove a slightly weaker version,
replacing the assertion that $H\in\bbb P_\tau$ with the assertion that
$H$ has infinitely many infinite vertical sections.  This weaker
version would be an immediate consequence of part~(2) if we knew that
there exists a sisnis ultrafilter, because any set in a sisnis
ultrafilter has infinitely many infinite sections.  The existence of a
sisnis ultrafilter, which is equivalent to the existence of infinitely
many non-isomorphic selective ultrafilters, is not provable in ZFC;
indeed, it is not provable that there exists even one selective
ultrafilter.  Nevertheless, we can still use part~(2) to obtain the
weakened part~(1) as follows.  Regardless of the existence or
non-existence of sisnis ultrafilters, we can pass to a forcing
extension of the universe in which the continuum hypothesis holds and
there are no new reals.  (It suffices to adjoin a generic subset of
$\omega_1$ with countable forcing conditions.) In the extension, there
are, thanks to the continuum hypothesis, plenty of selective
ultrafilters ($2^{2^{\aleph_0}}$ of them), so we have a sisnis
ultrafilter and therefore have the weakened part~(1) of the theorem.
But this result is a statement entirely about real numbers (note in
particular that the partition can be coded by a real number, as it
involves only an analytic set and its complement).  Since the forcing
extension didn't add reals, the same result holds in the original
universe, as required.

To pass from the weakened version of part~(1) to the original version
where $H$ is required to be in $\bbb P_\tau$, it suffices to recall
that every set with infinitely many infinite sections has a subset in
\bbb P and that $\bbb P_\tau$ is dense in \bbb P.  Therefore, we can
just replace the $H$ from the weakened version with a subset in $\bbb
P_\tau$ to complete the proof of part~(1).

Before proceeding to part~(3), we explain a technical strengthening of
part~(1) that will be used in the proof of part~(3).  Suppose we are
given, in addition to the $\omega$-type $\tau$ and the partition, a
subset $A$ of $\omega^2$ with infinitely many infinite sections.
Then, in the forcing extension used in the proof of part~(1), we can
choose the selective ultrafilters so that the resulting sisnis
ultrafilter \scr W contains $A$.  Then, when we apply part~(2) with this
sisnis ultrafilter, we can arrange for the homogeneous set $H$ to be a
subset of $A$; since both $A$ and $H$ are in \scr W, we can replace $H$
by its intersection with $A$.  The passage from the forcing extension
to the ground model and the shrinking of $H$ to put it into $\bbb
P_\tau$ preserve this arrangement.  Therefore, in part~(1) of the
theorem, we can get the homogeneous set to be included in any
prescribed $A$ that has infinitely many infinite sections.

Finally, we prove part~(3).  The preceding technical improvement of
part~(1) applies in particular to any $A\in\bbb P_\tau$.  So we have
that, for any partition into an analytic set and its complement, the
homogeneous sets $H$ are dense in $\bbb P_\tau$, so any $\bbb
P_\tau$-generic ultrafilter contains such a homogeneous set.  Recall
that $\bbb P_\tau$ is dense in \bbb P, so genericity is the same for
these two forcing notions, and we therefore have that every \bbb
P-generic ultrafilter contains a homogeneous set.  (We have tacitly
used the fact that forcing by \bbb P adds no new reals, so the pieces
of the partition, being analytic or coanalytic, are the same before
and after the forcing.)
\end{pf}

Theorem~\ref{analytic} is the analog, in our non-P-point context, of
part~(2) of Proposition~\ref{sel-ramsey} for selective ultrafilters.
We also have the following analog of part~(3) of that proposition.

\begin{thm} \label{hod} Suppose the universe is obtained from some
  ground model by L\'evy-collapsing to $\omega$ all cardinals below
  some Mahlo cardinal of the ground model.  Then the partition
  properties in Theorem~\ref{analytic} hold with $HOD\bbb R$ in place
  of analytic.
\end{thm}

\begin{pf}
We can proceed as in the proof of Theorem~\ref{analytic} with only the
following changes.  Instead of citing Theorem~7 of \cite{sel-homog},
we cite Corollary~11.1, which asserts (among other things) that the
L\'evy-Mahlo model satisfies Theorem~7 with $HOD\bbb R$ in place of
analytic.  Also, in the proof of part~(1), it is no longer necessary
to force to obtain the continuum hypothesis; the L\'evy-Mahlo model
satisfies the continuum hypothesis, so plenty of selective
ultrafilters are available in it.  
\end{pf}

Finally, we point out that part~(3) of Theorem~\ref{hod}, the part
about \bbb P-generic ultrafilters admits an easy converse, which could
be viewed as a sort of complete combinatorics.

\begin{prop} \label{comcom2} Suppose the universe is obtained from
  some ground model by L\'evy-collapsing to $\omega$ all cardinals
  below some Mahlo cardinal of the ground model.  Suppose further that
  \scr W is a non-P-point in standard position and that, for at least
  one $\omega$-type $\tau$, \scr W has the following partition
  property.  For any $HOD\bbb R$ partition of $\bbb P_\tau$ into two
  pieces, there is a set $H\in\scr W\cap\bbb P_\tau$ such that all of
  its subsets in $\bbb P_\tau$ lie in the same piece of the partition.
  Then \scr W is \bbb P-generic over $HOD\bbb R$.
\end{prop}

\begin{pf}
  Since $\bbb P_\tau$ is dense in \bbb P, it suffices to prove that
  \scr W intersects every dense $HOD\bbb R$ subset \scr D of $\bbb
  P_\tau$.  Because \scr D is dense, there cannot be any $H\in\bbb
  P_\tau$ (whether in \scr W or not) such that all its subsets in
  $\bbb P_\tau$ are outside \scr D.  By the assumed partition property
  of \scr W, we infer that there is $H\in\scr W\cap\bbb P_\tau$ such
  that $H$ lies in \scr D (and so do all its subsets in $\bbb P_\tau$,
  but we don't need this part of the result).  So $H$ witnesses that
  \scr W meets \scr D, as required. 
\end{pf}

Note that Proposition~\ref{comcom2} needs to assume the partition
property for only one $\omega$-type $\tau$.  Genericity follows, and
with genericity, the partition properties for all other $\omega$-types
also follow.

Proposition~\ref{comcom2} is a partial analog, in our non-P-point
context, of Corollary~\ref{comcom1} for selective ultrafilters.  A
more complete analog would result from a positive answer to the
following open problem.

\begin{qu}
  Suppose \scr W is a non-P-point in standard position, and suppose it
  is $(2,4)$-weakly Ramsey (and therefore $(n,T(n))$-weakly Ramsey for
  all $n$ by Corollary~\ref{alln}). Must it have the infinitary
  partition property in part~(2) of Theorem~\ref{analytic}?  If, in
  addition, the universe is a L\'evy-Mahlo model, must \scr W have the
  corresponding partition property for $HOD\bbb R$ partitions?
\end{qu}

\end{document}